\documentclass[
    a4paper,
    man,
    floatsintext]{glossaPX2}

\usepackage[T1]{fontenc}
\usepackage[american]{babel}
\usepackage[
  style=numeric,
  backend=biber,
  sorting=nyt
]{biblatex}
\NewBibliographyString{unpublished}
\DefineBibliographyStrings{english}{
  unpublished = {Unpublished},
}
\DefineBibliographyStrings{american}{
  unpublished = {Unpublished},
}

\usepackage[font={footnotesize,it}]{caption}
\usepackage{csquotes}
\usepackage{booktabs} 
\usepackage{tabularx}
\usepackage{linguex} 
\usepackage{cgloss}
\usepackage{amsmath}
\usepackage{amsthm}
\DeclareMathAlphabet{\mathcal}{OMS}{cmsy}{m}{n}
\usepackage{graphicx}
\graphicspath{ {./graphics/} }
\usepackage{ragged2e}
\usepackage{hyperref}
\usepackage[nameinlink,noabbrev]{cleveref}
\usepackage{enumitem}
\newcommand{\quotes}[1]{``#1''}
\usepackage[title]{appendix}
\usepackage{booktabs}
\usepackage{tikz}

\newtheorem{lem}{Lemma}
\newtheorem{thm}{Theorem}
\newtheorem{corr}{Corollary}
\newtheorem{prop}{Proposition}
\newtheorem{define}{Definition}

\newtheorem{conj}{Conjecture}

\crefname{thm}{theorem}{theorems}
\Crefname{thm}{Theorem}{Theorems}

\crefname{lem}{lemma}{lemmas}
\Crefname{lem}{Lemma}{Lemmas}

\crefname{corr}{corollary}{corollaries}
\Crefname{corr}{Corollary}{Corollaries}

\crefname{prop}{proposition}{propositions}
\Crefname{prop}{Proposition}{Propositions}

\crefname{define}{definition}{definitions}
\Crefname{define}{Definition}{Definitions}

\crefname{conj}{conjecture}{conjectures}
\Crefname{conj}{Conjecture}{Conjectures}

\crefname{obs}{observation}{observations}
\Crefname{obs}{Observation}{Observations}

\crefname{rmk}{remark}{remarks}
\Crefname{rmk}{Remark}{Remarks}

\newcommand{\m}{\mathcal}

\newcommand{\bi}{\overset{\leftrightarrow}}
\newcommand{\dc}{\mathrm{dc}}

\newcommand{\N}{\mathbb{N}}
\newcommand{\Z}{\mathbb{Z}}
\newcommand{\R}{\mathbb{R}}
\newcommand{\Q}{\mathbb{Q}}
\newcommand{\C}{\mathbb{C}}

\title[Properties of the Dichromatic Polynomial]{Properties of the Dichromatic Polynomial}
\author[Crew, Keller, Kuehr, Lin, Monroy Rodríguez, Williams] 
{\spauthor{Lila Crew\\
  \institute{University of Waterloo}}
  \AND
{\spauthor{Benjamin Keller\\
  \institute{Princeton University}
  }}
  \AND
\spauthor{Katharina Kuehr\\
  \institute{Minerva University}
  }
  \AND
\spauthor{Shao Yuan Lin\\
  \institute{University of Waterloo}
  }
  \AND
\spauthor{Megan Ixchel Monroy Rodríguez\\
  \institute{Universidad de Guanajuato}
  }
  \AND
\spauthor{Tabitha Williams\\
  \institute{The University of Texas at San Antonio}}
}

\begin{document}
\maketitle

\begin{abstract}
The chromatic polynomial is a well-studied function of graphs whose coefficients, factorization, and other properties are well-known to encode structural information about the underlying graph. An analogous notion of graph colouring exists for directed graphs which is also well-studied, but little is yet known about the corresponding dichromatic polynomial.

We investigate properties of the dichromatic polynomial, including its coefficient interpretations, root-dense regions, and asymptotic growth. We also show that chordal digraphs' polynomials split into linear factors in a way analogous to the undirected case, confirm a conjecture on dichromatic uniqueness posed by Gonz\'{a}lez-Moreno et. al., and study the density of dichromatic polynomials within certain larger families of polynomials.
\end{abstract}


\section{Introduction}
\label{intro}
Graph colourings have been a major area of study in the field of combinatorics since at least the early 1900s. In particular, Birkhoff \cite{birkhoff1912determinant} first showed that the number of ways to properly colour a given graph with $n$ colours is a polynomial in $n$, called the chromatic polynomial, which satisfies a recursive edge deletion-contraction property. In 1982, Neumann-Lara \cite{Neumann-Lara_1982} defined a generalization of graph colourings to digraphs: in the undirected case, a colouring is proper if there are no monochromatic edges, while in the directed case, it is proper if there are no monochromatic directed cycles. Under this interpretation, directed graph colouring is a generalization of undirected graph colouring, since an undirected graph may be represented as a directed graph with the same colouring rules by replacing each edge with a pair of arcs in opposite directions. In 2011, Harutyunyan \cite{Harutyunyan2011BrookstypeRF} established some interesting properties of the dichromatic polynomial, defined to be the Neumann-Lara analogue of the chromatic polynomial, and in 2022 Gonz\'{a}lez-Moreno et. al. \cite{DichromaticPoly} further investigated the dichromatic polynomial, including establishing a deletion-contraction relation.

 In this paper, we continue this line of research establishing the dichromatic polynomial as a graph-theoretic object of interest. After providing background (section 2), we use a relationship with chromatic polynomials to establish several interpretations for dichromatic coefficients (section 3). Previous work by Hochstättler et. al. \cite{Hochstaettler2021} gave explicit formulas for the first two coefficients and the last coefficient; our formula gives an interpretation for any coefficient in terms of chromatic polynomials of graphs formed by deleting and contracting arcs in a digraph, which yields many natural corollaries, including a simple explicit formula for the third dichromatic coefficient in terms of directed triangles, and a clear method to generalize this to further coefficients (although as with the chromatic polynomial, these formulas get more complex the further from the front the coefficient is). 
 
 Then, we consider the objects called \textit{perfect elimination digraphs} by Haskin and Rose \cite{haskins1973toward} and further studied by Meister and Telle \cite{chordalDigraphs}, which we call \textbf{chordal digraphs} (Section 4). As per their name, these digraphs have a perfect elimination ordering similar to that in chordal graphs. We show that additionally (and analogously to undirected chordal graphs) these digraphs have a dichromatic polynomial that splits linearly into factors determined by the degree sequence of the perfect elimination ordering. Furthermore we show that chordal digraphs are perfect, even using the strictest definition of digraph clique number as given by Aboulker et. al. \cite{aboulker2023clique}. 
 
 Next, we analyze the root-dense regions of the dichromatic polynomial (section 5). In particular, we conclude that unlike the case of the chromatic polynomial, which has real-root-free regions on the real line (famously including $(1, \frac{32}{27}]$ \cite{Jackson_1993}), the real and complex roots of dichromatic polynomials are dense in $\R$ and $\mathbb C$ respectively, even when restricted to oriented graphs (those with no two arcs between the same vertices in opposite directions). 
 
 We then resolve conjectures on dichromatic uniqueness of tournaments posed by Gonz\'{a}lez-Moreno et. al. \cite{DichromaticPoly} (section 6) before analyzing the density of dichromatic polynomials within families of integer-coefficient polynomials related to hyperplane arrangements (section 7). Finally, we conclude with a bound on its asymptotic growth as a function on the number of colours using an analogous argument of Seymour for the chromatic polynomial \cite{Seymour_1997a},  and provide further related conjectures of our own (section 8).

 \section{Background and Terminology}

 For brevity, we use the notation $[n] = \{1, \dots, n\}$.

 We assume standard terminology from (directed) graph theory, as given in e.g. \cite{bangdigraphs} and \cite{diestel}. In particular, a \textbf{graph} $G = (V,E)$ consists of a finite set $V$ of vertices, and a set $E$ of unordered pairs of vertices. A \textbf{directed graph} or \textbf{digraph} $D = (V,A)$ consists of a finite set $V$ of vertices, and a set $S$ of \textit{ordered} pairs of vertices. Given $C \subseteq V(D)$, the \textbf{subdigraph of $D$ induced by $C$}, denoted $D[C]$, is the digraph $C,A(C)$ where $A(C)$ consists of all edges of $A$ with both vertices in $C$. For brevity, when we refer to a cycle of a digraph, we mean a directed cycle unless otherwise specified. 

 A \textbf{digon} of a directed graph $D$ is a $2$-vertex cycle consisting of the arcs $(u,v)$ and $(v,u)$ for some $u,v \in V(D)$. An \textbf{asymmetric arc} $(u,v) \in A(D)$ is an arc such that $(v,u) \notin A(D)$; equivalently, it is an arc that is not part of a digon.

 A \textbf{$k$-colouring of a graph $G$} is a map $c: V(G) \rightarrow \{1, \dots, k\}$. A $k$-colouring of a graph is \textbf{proper} if whenever $uv \in E(G)$ we have $c(u) \neq c(v)$. The \textbf{chromatic polynomial} $\chi_G(x)$ of $G$ is the unique function such that for each positive integer $k$, we have $\chi_G(k)$ is the number of proper $k$-colourings of $G$.

 A \textbf{$k$-colouring of a directed graph (or digraph) $D$} is a map $d: V(D) \rightarrow \{1, \dots, k\}$. A $k$-colouring of a directed graph is \textbf{proper} if $d$ does not assign the same colour to every vertex of any directed cycle of $D$. The \textbf{dichromatic polynomial} $P_D(x)$ of $D$ is the unique function such that for each positive integer $k$, we have $P_D(k)$ is the number of proper $k$-colourings of $D$ (that this is a polynomial function is proved in e.g. \cite{DichromaticPoly}). 

\section{Coefficients of the Dichromatic Polynomial}

Determining a combinatorial interpretation of the coefficients of the chromatic polynomial is a classical problem in graph theory. For instance, if $G$ is a graph with $n$ vertices, $m$ edges, $t$ triangles, and $\chi_G(x)=\sum_{i=0}^{n}a_ix^i$, then $a_n=1,a_{n-1}=-m,a_{n-2}=\binom{m}{2}-t$, and $a_0=0$ \cite{Whitney_1932}. Specifically, the first few chromatic coefficients encode specific structural information about $G$. In general, combinatorial interpretations for arbitrary chromatic coefficients are given via Whitney's Broken Circuit Theorem (\cite{Whitney_1932}).

Initial results of this kind are also known for the dichromatic polynomial. In \cite{DichromaticPoly}, it was shown that if $D$ is a digraph of order $n$, $q$ is the number of digons in $D$, and $P_D(x)=\sum_{i=0}^{n}c_ix^i$, then $c_n=1$, $c_{n-1}=-q$, and $c_0=0$. In particular, the leading and constant coefficients behave as in the chromatic setting, while digons play the role of edges in the second coefficient.

Building on these results, we study the coefficients of the dichromatic polynomial in greater detail. Unlike in ordinary graph colouring, an individual asymmetric arc imposes no restriction on a colouring; rather, the fundamental obstructions are the monochromatic directed cycles. As a result, the coefficients of the dichromatic polynomial $P_D(x)$ are governed by the minimal vertex sets which induce directed cycles and by the ways in which these sets intersect. Digons behave like ordinary edges, while directed cycles involving asymmetric arcs must be accounted for differently. By separating these two types of obstructions, we obtain a formula for $P_D(x)$ in terms of chromatic polynomials of contractions of the symmetric-pair graph $S(D)$ (defined below). We then use this formula to determine when the leading chromatic and dichromatic coefficients agree, identify the first point of deviation, and derive an explicit formulation for the third dichromatic coefficient.

\vspace{\baselineskip}

To formalize the cyclic obstructions and their intersections, we begin with the following definitions:

\begin{define}\label{define:min_cyclic_sets}
Let $D=(V,A)$ be a digraph. Then let \[\m C(D)=\{C\subseteq V \mid D[C]\text{ contains a cycle}, \nexists C'\subset C\text{ satisfying }D[C']\text{ contains a cycle}\}\]

which is the set of minimal cyclic vertex sets of $D$.
\end{define}

\begin{define}\label{define:equivalence_relation}
Let $D=(V,A)$ be a digraph and let $\m F\subseteq\m C(D)$. Define the equivalence relation $\sim_\m F$ on $V$ as follows: $u\sim_\m F v$ if and only if either $u=v$ or there exists a chain $C_1,\dots, C_m$ such that \[u\in C_1, C_1\cap C_2\neq\emptyset, C_2\cap C_3\neq\emptyset,\dots,v\in C_m\]

where each $C_i\in\m F$. In other words, $u\sim_\m F v$ if $u,v$ are forced to have the same colour if we force all directed cycles in $\m F$ to be monochromatic.
\end{define}

\begin{define}\label{define:equivalence_classes}
Let $D=(V,A)$ be a digraph and let $\m F\subseteq\m C(D)$. Let $\pi(\m F)$ denote the partition of $V$ into equivalence classes given by $\sim_\m F$ and let $\kappa(\m F)=|\pi(\m F)|$ be the number of equivalence classes.
\end{define}

\begin{define}\label{define:symmetric_pair_graph}
Let $D=(V,A)$ be a digraph. Then the \textnormal{symmetric-pair digraph}, denoted by $S(D)$, has vertex set $V$ and edges \[E=\{uv\mid (u,v),(v,u)\in A\}.\]
\end{define}

With these definitions in mind, we prove the aforementioned identity relating the dichromatic polynomial to chromatic polynomials of contractions of the symmetric-pair digraph $S(D)$.

\begin{thm}\label{thm:dichromatic_chromatic_identity}
Let $D=(V,A)$ be a digraph. Then \[P_D(x)=\sum_{\m Y\subseteq\m H}(-1)^{|\m Y|}\chi_{S(D)/\pi(\m Y)}(x)\]

where $\m H=\{C\in\m C(D)\mid |C|\geq 3\}$, $\pi(\m Y)$ is the partition of $V$ obtained by forcing every $C\in\m Y$ to be monochromatic via the equivalence class $\sim_{\m Y}$, and $S(D)/\pi(\m Y)$ is the graph obtained from $S(D)$ by contracting every block of $\pi(\m Y)$ into a single vertex. (For convention, if an edge of $S(D)$ has endpoints within the same block, contracting it creates a loop, and we can allow loops with the convention that a graph with loops has chromatic polynomial $0$).
\end{thm}

\begin{proof}
Recall from \cite{Harutyunyan2011BrookstypeRF} (Section 8.4) and \cite{Dohmen_Trinks_2014} that \[P_D(x)=\sum_{\m F\subseteq\m C(D)}(-1)^{|\m F|}x^{\kappa(\m F)}\]

and observe that every $\m F\subseteq\m C(D)$ decomposes uniquely as $\m F=\m E\cup\m Y$ where $\m E$ is a set of digons (corresponding to some set $\m E'\subseteq E(S(D))$) and $\m Y\subseteq\m H$ is a family of digon-free minimal cyclic sets. This gives us that \[P_D(x)=\sum_{\m Y\subseteq\m H}(-1)^{|\m Y|}\sum_{\m E'\subseteq E(S(D))}(-1)^{|\m E'|}x^{\kappa(\m E\cup\m Y)}\]

Now fix $\m Y\subseteq\m H$ and let $G_\m Y=S(D)/\pi(\m Y)$ (notice that $G_\m Y$ may be a multigraph with loops). Clearly, the vertices of $G_\m Y$ are the equivalence classes (blocks) of $\pi(\m Y)$. Furthermore, each digon in $D$ (equivalently, each edge in $S(D)$) maps to an edge of $G_\m Y$, possibly a loop, so for $\m E'\subseteq E(S(D))$, let $\overline{\m E'}$ be its image in $G_\m Y$. Then imposing the \quotes{bad events} in $\m E\cup\m Y$ gives exactly the same equivalence relation as first contracting the blocks of $\pi(\m Y)$ and then contracting the edges $\overline{\m E'}$. Hence, \[\kappa(\m E\cup\m Y)=c((V(G_\m Y),\overline{\m E'})\]

This gives us that \[\sum_{\m E'\subseteq E(S(D))}(-1)^{|\m E'|}x^{\kappa(\m E\cup\m Y)}=\sum_{\overline{\m E'}\subseteq E(G_\m Y)}(-1)^{|\overline{\m E'}|}x^{c((V(G_\m Y),\overline{\m E'})}=\chi_{G_\m Y}(x)\]

So, we get that \[P_D(x)=\sum_{\m Y\subseteq\m H}(-1)^{|\m Y|}\chi_{G_\m Y}(x)=\sum_{\m Y\subseteq\m H}(-1)^{|\m Y|}\chi_{S(D)/\pi(\m Y)}(x)\]

as desired.
\end{proof}

\Cref{thm:dichromatic_chromatic_identity} immediately gives us several corollaries, including polynomial identities and an explicit formula for the third dichromatic coefficient. 

\begin{corr}\label{corr:coefficient_identity}
Let $D=(V,A)$ be a digraph of order $n$ and let $j\in[n]$. Then \[[x^{n-j}]P_D(x)=[x^{n-j}]\chi_{S(D)}(x)+\sum_{\emptyset\neq\m Y\subseteq\m H}(-1)^{|\m Y|}[x^{n-j}]\chi_{S(D)/\pi(\m Y)}(x).\]
\end{corr}

\begin{corr}[Alignment Theorem]\label{corr:alignment}
Let $D=(V,A)$ be a digraph of order $n$ and suppose that every directed cycle of $D$ using at least one asymmetric arc has length at least $g$. Then for $0\leq j\leq g-2$, we have that \[[x^{n-j}]P_D(x)=[x^{n-j}]\chi_{S(D)}(x)\]

where $\chi_{S(D)}(x)$ is the ordinary chromatic polynomial for the undirected graph $S(D)$.
\end{corr}

\begin{corr}[First Deviation Theorem]\label{corr:first_deviation}
Let $D=(V,A)$ be a digraph of order $n$ and suppose that every directed cycle of $D$ using at least one asymmetric arc has length at least $g$. Then \[[x^{n-g+1}]P_D(x)=[x^{n-g+1}]\chi_{S(D)}(x)-N_g(D)\]

where $N_g(D)$ counts the number of directed $g$-cycles whose vertices do not induce any digons.
\end{corr}

\begin{corr}[Third Dichromatic Coefficient]\label{corr:third_coefficient}
Let $D=(V,A)$ be a digraph of order $n$ with dichromatic polynomial $P_D(x)=c_nx^n+c_{n-1}x^{n-1}+\cdots+c_1x+c_0$. Then \[c_{n-2}=\binom{q}{2}-t_S-t_D\]

where $q$ is the number of digons, $t_S$ is the number of triangles in $S(D)$, and $t_D$ is the number of directed triangles whose vertices do not induce a digon.
\end{corr}

\section{Chordal Digraphs}

The class of chordal graphs, most directly defined as graphs in which there are no induced cycles of length larger than $3$ (so every large cycle has a chord shortcutting it) is foundational in graph theory. These graphs admit a number of nice properties; however, it is not immediately obvious what a correct analogue of the family is in digraphs. One obvious choice is to take chordal digraphs to be those in which every directed cycle has a chord; however, such digraphs fail to have many of the nice properties of chordal graphs. For example, they are not perfect, and in fact they are not even $\chi$-bounded with respect to any accepted definition of digraph clique number \cite{aboulker2022chordal}.

In this section we further explore a class of digraphs, the \textbf{perfect elimination digraphs} of Haskin and Rose \cite{haskins1973toward}, that have previously been considered as a natural candidate for a digraph analogue of chordal graphs (e.g. by Meister and Telle \cite{chordalDigraphs}). In particular, we show that this class of digraphs plays nicely with the dichromatic polynomial in a way that is analogous to the case of undirected chordal graphs. Hence, throughout this section we will simply refer to them as \textbf{chordal digraphs}.

First, we define these chordal digraphs. One characterization of undirected chordal graphs is that they admit a \textit{perfect elimination ordering}, in the sense that they can be reduced to a single vertex by repeatedly deleting \textit{simplicial vertices} in the remaining graph. Thus, we first give a notion of directed simplicial vertex, which requires an appropriate notion of a clique in a directed graph.

\begin{define}[\cite{chordalDigraphs}]\label{define:chordal_digraphs}

\hspace{-1.5em}

    \begin{enumerate}
    
        \item A \textbf{di-clique of a digraph $D$} is a pair $(X,Y)$ with $X,Y\subset V(D)$ where for all $x\in X$ and $y\in Y$ with $x\neq y$, $(x,y)\in A(D)$.
        \item The class of \textbf{chordal digraphs} is inductively defined as follows:
        \begin{enumerate}
            \item A digraph on a single vertex is a chordal digraph.
            \item Let $D$ be a chordal digraph and let $u$ be a vertex that does not appear in $D$. Let $(X,Y)$ be a di-clique of $D$. The digraph obtained from $D$ by adding $u$ and the set of arcs $\{(x,u):x\in X\}\cup\{(u,y):y\in Y\}$ is also a chordal digraph.
        \end{enumerate}
    \end{enumerate}
    
\end{define}

\begin{define}[\cite{chordalDigraphs}]\label{define:di_simplicial}
    A vertex $u$ of a digraph $D$ is \textbf{di-simplicial} if $(N_D^{\mathrm{in}}(u),N_D^{\mathrm{out}}(u))$ is a di-clique of $D$.
\end{define}

That is, a chordal digraph is built up by repeatedly adding di-simplicial vertices. Thus, chordal digraphs have a perfect elimination order, formed by deleting di-simplicial vertices in the reverse of the order they were added. In fact, as in chordal graphs, you don't even need to know the original order to do this.


\begin{thm}\cite{chordalDigraphs}\label{thm:peo}
    A digraph $D$ is chordal if and only if $D$ can be reduced to a digraph on a single vertex by repeatedly deleting an \textnormal{arbitrary} di-simplicial vertex.
\end{thm}

We can now understand dichromatic polynomials of chordal digraphs by determining how they change upon deleting a di-simplicial vertex.

\begin{prop}\label{prop:chordal_decomposition}
Let $v$ be a di-simplicial vertex of a digraph $D$ and $Q:=N_D^{\mathrm{in}}(v)\cap N_D^{\mathrm{out}}(v)$. Then \[P_D(k)=(k-q)P_{D-v}(k),\]

where $q=|Q|$.
\end{prop}

\begin{proof}
Let $C$ be a valid $k$-colouring in $D-v$. Notice that if adding $v$ and colouring it arbitrarily creates a monochromatic cycle, then we have two cases.

\begin{description}
    \item[Case 1:] If the cycle has length greater than or equal to $3$, then it is of the form $v\to u\to\cdots\to w\to v$. However, since $v$ is di-simplicial, we also have that $(w,u)\in A(D)$, meaning that we have the monochromatic directed cycle $ u\to\cdots\to w\to u$, which is a monochromatic directed cycle in $D-v$. Since we assumed that $C$ was a valid colouring of $D-v$, this is a contradiction.
    \item[Case 2:] If the cycle has length $2$ (i.e., if $v$ forms a digon with some vertex $u$), then it must be the case that $u$ and $v$ have different colours. Now, consider distinct vertices $u_1,u_2\in Q$. Since $v$ is di-simplicial, we have that $(u_1,u_2),(u_2,u_1)\in A(D)$, meaning that $u_1$ and $u_2$ are connected with a digon and therefore must also have different colours in the original colouring of $D-v$. Thus all vertices of $Q$ have different colours, implying that there are exactly $q$ forbidden colours for $v$.

    \vspace{5pt}

    Furthermore, any vertex in $(N_D^{\mathrm{in}}(v)\cup N_D^{\mathrm{out}})\setminus Q$ cannot form a digon with $v$. Since Case 1 implies that monochromatic cycles of length $\geq 3$ cannot exist, the colours of vertices not in $Q$ impose no colouring constraints on $v$. Thus, there are exactly $q$ forbidden colours for $v$, leaving exactly $k-q$ choices for colouring it, indepedent of the original colouring of $D-v$. Hence, \[P_D(k)=(k-q)P_{D-v}(k).\]
\end{description}

\end{proof}

Since chordal digraphs have a perfect elimination order, we can repeatedly apply the above proposition until the dichromatic polynomial splits into linear factors, analogously to the undirected case.

\begin{corr}\label{corr:chordal_dichromatic_factors}
    Let $D$ be a chordal digraph with perfect elimination order $v_1,\dots, v_n$ (meaning each $v_i$ is a di-simplicial vertex of $D-\{v_1, \dots, v_i\}$). Then \[P_D(k)=\prod_{i\in [n]}(k-q_i)\]

    where $q_i=|Q_i|$ for \[Q_i:=N_{D_i}^{\mathrm{in}}(v_i)\cap N_{D_i}^{\mathrm{out}}(v_i),\]
    
    and $D_i$ is the subdigraph induced by $v_i, \dots, v_n$.
\end{corr}

\medskip

Thus, dichromatic polynomials of these digraphs split linearly into factors determined by a perfect elimination ordering. Now, we can answer questions about the dichromatic number of chordal digraphs. 

In particular, every undirected chordal graph $G$ is perfect, meaning that its chromatic number $\chi(G)$ (the minimum $k$ such that $G$ has a $k$-colouring) is equal to its clique number $\omega(G)$ (the size of the largest set $A \subseteq V(G)$ such that every pair of vertices in $A$ is an edge of $G$). It is straightforward to see that for every graph $G$ we have $\omega(G) \leq \chi(G)$, so perfect graphs are those where the inequality is an equality. Our goal in the remainder of this section is to show that chordal digraphs are perfect, meaning their dichromatic number equals their clique number. We will need to give appropriate definitions of each of these.

The \textbf{dichromatic number of a digraph $D$}, notated $\dc(D)$, is the smallest positive integer $k$ such that $D$ admits a $k$-colouring; equivalently, it is the smallest positive integer $k$ such that $P_D(k) > 0$. Then the following is immediate.

\begin{corr}\label{corr:dichromatic_number}
    Let $D$ be a chordal digraph with roots $q_i$ as in Corollary \ref{corr:chordal_dichromatic_factors} for $i\in[n]$. Then $\dc(D)=\max(q_i)+1$.
\end{corr}


As far as clique number of a digraph goes, there are a few distinct definitions. The definition taken by Andres and H\"{o}chstattler \cite{perfectDigraphs} is a common one.

\begin{define}[\cite{perfectDigraphs}]\label{define:largest_digon_vertices}
    Let $D$ be a digraph. Then the Andres-Hochstättler clique number $\omega_{AH}(D)$ is the largest number of vertices of a subset $A$ of $V(D)$ such that every pair of vertices in $A$ induces a digon.
\end{define}

Equivalently, this is the size of the largest clique in $S(D)$.

 Another definition of clique number in directed graphs is due to Aboulker et. al. \cite{aboulker2023clique} and arises from the notion of (undirected) backedge graphs of a digraph.

\begin{define}\label{define:backedge_graph}
    Let $D=(V,A)$ be a digraph and $\pi=(v_1,\dots,v_n)$ be an ordering of $V$. We define the (undirected) backedge graph $B_\pi=(V,E)$ where $\{v_i,v_j\}\in E$ if $(v_j,v_i)\in A$ for $j>i$. 
\end{define}

\begin{define}[\cite{aboulker2023clique}]\label{define:clique_num}
    Let $D$ be a digraph. We define the backedge clique number $\omega_B(D)$ of $D$ by \[\omega_B(D)=\min_\pi \omega(B_\pi).\]
\end{define}

Certainly a clique in $S(D)$ will form a clique in every backedge graph of $D$, so it is clear that $\omega_{AH}(D) \leq \omega_{B}(D)$ (and the equality can be strict, since for example tournaments $T$ are oriented complete graphs so have $\omega_{AH}(T) = 1$, but can have arbitrarily large backedge clique number \cite{aboulker2023clique}). It is likewise straightforward to show that $dc(D) = \min_{\pi} \chi(B_{\pi})$ (e.g. \cite{aboulker2023clique}), so it is clear that $\omega_B(D) \leq dc(D)$. Thus, since $\omega_{AH}(D) \leq \omega_B(D) \leq dc(D)$, to show that chordal digraphs are perfect with respect to both notions of clique number, it is enough to show that they are perfect with respect to the smaller one.

\begin{prop}\label{prop:dichromatic_num_equals_digon_num}
    Let $D$ be a chordal digraph. Then \[\dc(D)=\omega_{AH}(D).\]
\end{prop}

\begin{proof}
As above, $\omega_{AH}(D) \leq\dc(D) = \max(q_i)+1$.

Now, let $v_1, \dots, v_n$ be a perfect elimination order of $D$, and let $u$ be a vertex with $q_u=\max(q_i)$. Since $u$ is a di-simplicial vertex, every pair of vertices in $Q_u$ forms a digon. Thus, there is a clique $C'$ in $D$ such that every pair of vertices in $C'$ induces a digon with $q_u+1$ vertices. Hence $\omega_{AH}(D)\geq q_u+1$. Combining this with the previous inequality yields that $\omega_{AH}(D)=\max(q_i)+1$. So, \Cref{corr:dichromatic_number} implies that $\omega_{AH}(D)=\max(q_i)+1=\dc(D)$, as desired.
\end{proof}






\section{Roots of the Dichromatic Polynomial}
The location and distribution of roots of the chromatic polynomial have long been central topics in algebraic graph theory. Although complex chromatic roots are dense in $\C$, real chromatic roots are much more restricted (\cite{SOKAL_2004}). In particular, no chromatic roots exist in the intervals $(-\infty,0)$, $(0,1)$, or $(1,32/27]$, and real chromatic roots are dense in all other regions of $\R$ (\cite{Jackson_1993},\cite{THOMASSEN_1997}). This interesting behavior raises the analogous question of how roots of the dichromatic polynomial behave. In this section, we show that dichromatic roots exhibit substantially greater flexibility. Our main results establish that real dichromatic roots are dense in $\R$ and that complex dichromatic roots are dense in $\C$. Moreover, both density results continue to hold when restricted to digon-free digraphs. 

\vspace{\baselineskip}

Our approach first relates the dichromatic polynomials of a particular family of digraphs to hypergraph chromatic polynomials and, in turn, to independence polynomials of bipartite graphs. We then introduce constructions that shift real dichromatic roots while preserving the relevant structural conditions. 

Let $G=(V,E)$ be a simple graph. A \textbf{hypergraph} $H = (V,E)$ consists of a vertex set $V$ and a finite set $E$ of unordered subsets of $V$ called hyperedges. A \textbf{$k$-colouring of a hypergraph $H$} is a map $f: V(H) \rightarrow [k]$ such that no hyperedge is monochromatic. That is, there is no $e \in E$ such that $f(v)$ is the same for all $v \in e$. Recall from \cite{Zhang_Dong_2017} the definition of the hypergraph $\m H_{\bullet G}$, whose vertices are $V(G)\cup \{w\}$ and whose hyperedges are $\{u,v,w\}$ for each edge $e=uv$. 

\begin{lem}\label{lem:bipartite_digraph_hypergraph_lemma}
Let $G=(V,E)$ be a simple bipartite graph with bipartition $(A,B)$ and define a digraph $D_G$ on $V(G)\cup\{w\}$ as follows: for every edge $ab\in E(G)$ with $a\in A,b\in B$, add the arcs $w\to a,a\to b,b\to w$ (so every edge $ab$ gives the directed triangle $w\to a\to b\to w$). Then, \[P_{D_G}(x)=P_{\m H_{\bullet G}}(x)\]
\end{lem}

\begin{proof}
If $\{a,b,w\}$ is a monochromatic hyperedge in $\m H_{\bullet G}$, then since $ab\in E(G)$, the construction of $D_G$ induces the directed triangle $w\to a\to b\to w$, which is a monochromatic directed cycle in $D_G$ since $a,b,w$ have the same colour (they are vertices of a monochromatic hyperedge). Conversely, if $C$ is a monochromatic directed cycle in $D_G$, then $C$ must follow the pattern: \[w\to A\to B\to w\]

since the only arcs are of the form $w\to A$, $A\to B$, and $B\to w$. As such we can see that $C$ must be a triangle of the form $w\to a\to b\to w$ for some edge $ab\in E(G)$. Since $C$ is monochromatic, $a,b,w$ have the same colour, so this creates the monochromatic hyperedge $\{a,b,w\}$ in $\m H_{\bullet G}$.

\vspace{10pt}

This establishes a bijection between the improper colourings of $D_G$ and $\m H_{\bullet G}$ (a colouring $f$ of $V(G)\cup\{w\}$ has a monochromatic directed cycle in $D_G$ if and only if $f$ creates a monochromatic hyperedge in $\m H_{\bullet G}$), which implies a bijection between the proper colourings of the two graphs. This gives us that \[P_{D_G}(x)=P_{\m H_{\bullet G}}(x).\]
\end{proof}

We now define the independence polynomial and the lexicographic product (or composition) of two graphs.

\begin{define}\label{define:independent_set/polynomial}
An independent set in a graph $G$ is a set of vertices which do not share any edges among each other. The independence number $\beta$ of $G$ is the size of the largest independent set of $G$. The independence polynomial of $G$ is denoted by $I(G,x)=\sum_{k=0}^{\beta}i_kx^k$, where $i_k$ is the number of independent sets of size $k$.
\end{define}

\begin{define}\label{define:lexicographic_product}
For two graphs $G$ and $H$, let $G\times H$ be the graph with vertex set $V(G)\times V(H)$ and such that vertex $(a,x)$ is adjacent to vertex $(b,y)$ if and only if a is adjacent to $b$ (in $G$) or $a=b$ and $x$ is adjacent to $y$ (in $H$). The graph $G\times H$ is the lexicographic product (or composition) of $G$ and $H$, and can be thought of as the graph arising from $G$ and $H$ by substituting a copy of H for every vertex of $G$.
\end{define}

\begin{lem}\label{lem:independence_root_density}
Independence roots of simple bipartite graphs are dense in $(-\infty,0]$.
\end{lem}

\begin{proof}
Theorem 8 in \cite{Brown_Hickman_Nowakowski_2004} tells us that independence roots of paths $P_n$ (which are clearly simple and bipartite) are dense in $(-\infty,-1/4]$, so all that remains to show is that independence roots of simple bipartite graphs are dense in $(-1/4,0)$.

\vspace{10pt}

The paper \cite{Brown_Hickman_Nowakowski_2004} gives the following identity, where $I(G,x)$ is the independence polynomial of the graph $G$: \[I(G\times H,x)=I(G,I(H,x)-1).\]

Taking $H=\overline{K_m}$ (the complement of the complete graph on $m$ vertices), we have that \[I(\overline{K_m},x)=\sum_{k=0}^{m}\binom{m}{k}x^k=(1+x)^m\]

since each set of vertices is independent. Choosing $G=P_n$ (the path of length $n$), we get that \[I(P_n\times \overline{K_m},x)=I(P_n,I(\overline{K_m},x)-1)=I(P_n,(1+x)^m-1)\]

Also notice that $P_n\times \overline{K_m}$ is bipartite since it is simply the path $P_n$ with each vertex replaced by an independent set of size $m$, and each edge replaced by a complete bipartite graph between the corresponding independent sets. Now let $s\in(-1/4,0)$, $\epsilon>0$, and choose an $m\in\N$ sufficiently large so that \[(1+s)^m-1\in(-1,-1/4)\]

This is possible since $1+s\in(3/4,1)$, so $(1+s)^m-1\to -1$ as $m\to\infty$. Set $y=(1+s)^m-1$. By Theorem 8 of \cite{Brown_Hickman_Nowakowski_2004}, path independence roots are dense in $(-\infty,-1/4]$, so choose a path root $r$ of some $I(P_n,x)$ very close to $y$ (in particular, choose it close enough such that $r\in (-1,-1/4)$ and $(1+r)^{1/m}-1$ is within $\epsilon$ of $s$). This is possible because the map $r\mapsto (1+r)^{1/m}-1$ is continuous on $(-1,0)$. Now define $x=(1+r)^{1/m}-1$. Then, $r=(1+x)^m-1$. So, \[I(P_n\times \overline{K_m},x)=I(P_n,(1+x)^m-1)=I(P_n,r)=0\]

and $x$ is within $\epsilon$ of $s$. Thus, the roots of the bipartite graphs $P_n\times \overline{K_n}$ are dense in $(-1/4,0)$. Combining this result with Theorem 8 yields the lemma.
\end{proof}

\begin{corr}\label{corr:real_dichromatic_partial_density}
Real dichromatic roots are dense in $(-\infty,1]$.
\end{corr}

\begin{proof}
Let $G=(V,E)$ be a simple bipartite graph with bipartition $(A,B)$ and let $D_G$ be as in \Cref{lem:bipartite_digraph_hypergraph_lemma}. Theorem 1 in \cite{Zhang_Dong_2017} gives us that \[P_{\m H_{\bullet G}}(x)=x(x-1)^{|V(G)|}I\left(G,\frac{1}{x-1}\right)\]

Now let $\alpha<0$ be an independence root of $G$; so $I(G,\alpha)=0$. Set $\alpha=\frac{1}{x-1}\implies x=1+\frac{1}{\alpha}$. Then \Cref{lem:bipartite_digraph_hypergraph_lemma} implies that \[P_{D_G}(x)=P_{\m H_{\bullet G}}(x)=x(x-1)^{|V(G)|}I\left(G,\frac{1}{x-1}\right)=x(x-1)^{|V(G)|}I(G,\alpha)=0\]

So, every negative independence root $\alpha$ resulting from a simple bipartite graph produces a dichromatic root $x=1+\frac{1}{\alpha}$. Since the $\alpha$'s are dense in $(-\infty,0]$ (\Cref{lem:independence_root_density}), the corresponding $x$'s are dense in $(-\infty,1]$.
\end{proof}

\begin{corr}\label{corr:real_dichromatic_entire_density}
Real dichromatic roots are dense in $\R$.
\end{corr}

\begin{proof}
Given a digraph $D$, construct the digraph $J_k(D)$ by adding $k$ new vertices and putting digons between every pair of new vertices and between every new and every old vertex. Then, \[P_{J_k(D)}(x)=x^{\underline{k}}P_D(x-k)\]

where $x^{\underline k}=x(x-1)\cdot\cdots\cdot(x-k+1)$ is the $k$-th falling factorial. This is the case since each of the $k$ new vertices must receive $k$ distinct colours (due to the digon constraints), and all old vertices must avoid these $k$ colours (again due to the digon constraints). So, the old copy of $D$ is coloured acyclically using the remaining $x-k$ colours. Thus, if $\rho$ is a root of $P_D(x)$, then $\rho+k$ is a root of $P_{J_k(D)}(x)$. Since the original roots $\rho$ are dense in $(-\infty,1)$ (\Cref{corr:real_dichromatic_partial_density}), their shifts $\rho+k$ for $k\in\{0,1,2,\dots\}$ are dense in \[\bigcup_{k\in\{0,1,2,\dots\}}(-\infty,1+k)=\R.\]

So, real dichromatic roots are dense in $\R$.
\end{proof}

\begin{corr}\label{corr:complex_dichromatic_density}
Complex dichromatic roots are dense in $\C$. Furthermore, complex dichromatic roots of digon-free digraphs are dense in $\C$.
\end{corr}

\begin{proof}
The density of complex dichromatic roots follows from the density of complex chromatic roots and the fact that the chromatic and dichromatic polynomials are equal for a graph and its biorientation, respectively. The second statement follows from the same construction as was used in the density of real dichromatic roots.
\end{proof}

It remains to show that real dichromatic roots of digon-free digraphs are dense in $\R$. We first need the following auxiliary lemma.

\begin{lem}\label{lem:close_sign_changing_root}
Let $D=(V,A)$ be a digon-free digraph of order $n$ and let $\alpha>0$ be a real root of $P_D(x)$ at which $P_D$ changes sign. Then, for every $\epsilon>0$, there exists a digon-free digraph $D'$ whose dichromatic polynomial has a sign-changing real root $\beta$ which satisfies \[|\beta-(\alpha+1)|<\epsilon\]
\end{lem}

\begin{proof}
For $t\geq 1$, define $J_t(D)$ as follows: retain $D$, add a new vertex $u$, and for every vertex $v\in V$, add the vertices $z_{v,1},\dots,z_{v,t}$ and add the arcs $u\to v$, $v\to z_{v,j}$, and $z_{v,j}\to u$. Thus, every $u,v,z_{v,j}$ spans a directed triangle. Furthermore, since every pair receives arcs in at most one direction, it is clear that $J_t(D)$ contains no digons.

\vspace{10pt}

Now define $\m A(D)=\{S\subseteq V\mid D[S]\text{ is acyclic}\}$. We now claim that \begin{equation}P_{J_t(D)}(q)=q\sum_{S\in\m A(D)}P_{D-S}(q-1)\cdot (q-1)^{t|S|}q^{t(n-|S|)}\end{equation}

To prove this, let $q$ be a positive integer and choose the colour of some vertex $u$ (note that there are $q$ ways to do this). Then define $S=\{v\in V\mid v\text{ has the same colour as }u\}$, which is necessarily an acyclic vertex set. Thus, the vertices in $D-S$ may be coloured acyclically with the other $q-1$ colours in $P_{D-S}(q-1)$ ways. Now, for each $v\in S$, each $z_{v,j}$ has $q-1$ choices since it cannot receive the common colour of $u$ and $v$ (since then $u\to v\to z_{v,j}\to u$ would be a monochromatic cycle). Finally, for $v\notin S,$, $z_{v,j}$ has all $q$ choices.

These are the only restrictions involving the new vertices. Indeed, any directed cycle containing $z_{v,j}$ contains $v\to z_{v,j}\to u$, so if such a cycle were monochromatic, then $u,v,z_{v,j}$ would have the same colour, and the directed triangle on these three vertices would already be monochromatic. This proves (1).

Dividing this polynomial identity by $q^{tn}$ yields \begin{equation}q^{-tn}P_{J_t(D)}(q)=q\sum_{S\in\m A(D)}P_{D-S}(q-1)\left(\frac{q-1}{q}\right)^{t|S|}\end{equation}

Now, notice that on every compact interval $K\subset(1,\infty)$, \[\max_{q\in K}\left|\frac{q-1}{q}\right|<1\]

Hence, all terms of (2) with $S\neq\emptyset$ converge uniformly to $0$ as $t\to\infty$. The $S=\emptyset$ term is just $qP_D(q-1)$, so \begin{equation}q^{-tn}P_{J_t(D)}(q)\overset{t\to\infty}{\longrightarrow} qP_D(q-1)\end{equation}

uniformly on compact subintervals of $(1,\infty)$. Now choose $0<a<\alpha<b$, arbitrarily close to $\alpha$, such that $P_D(a)P_D(b)<0$. Then, $(1+a)P_D(a)$ and $(1+b)P_D(b)$ have opposite signs. Furthermore, by (3), for all sufficiently large $t$, $P_{J_t(D)}(1+a)$ and $P_{J_t(D)}(1+b)$ have opposite signs. Hence, $P_{J_t(D)}$ has a real root $\beta\in(1+a,1+b)$, and since the signs at the endpoints differ, at least one root in this interval has odd multiplicity and is therefore sign-changing. Thus, taking $a,b$ sufficiently close to $\alpha$ gives $|\beta-(\alpha+1)|<\epsilon$, as desired.
\end{proof}

\begin{thm}\label{thm:real_digon_free_density}
Real dichromatic roots of digon-free digraphs are dense in $\R$.
\end{thm}

\begin{proof}
Fix $x\in\R$ and $\epsilon>0$. If $x\leq1$, then \Cref{corr:real_dichromatic_partial_density} gives a digon-free digraph with a real dichromatic root within $\epsilon$ of $x$. Now suppose $x>1$. Choose an integer $k\geq 0$ and $y\in (0,1)$ such that \[|x-(y+k)|<\frac{\epsilon}{2}\]

If $x\notin\Z$, then take $k=\lfloor x\rfloor$ and $y=x-k$. If $x\in\Z$, then take $k=x-1$ and choose $y<1$ sufficiently close to $1$. Now choose \[0<\eta<\min\left\{\frac{y}{2},\frac{\epsilon}{2(k+1)}\right\}\]

Then by \Cref{corr:real_dichromatic_partial_density}, there exists a digon-free digraph $D_0$ with a real root $\alpha_0\in (0,1)$ satisfying $|\alpha_0-y|<\eta$. Now apply \Cref{lem:close_sign_changing_root} repeatedly: for $i=1,\dots,k$, construct a digon-free digraph $D_i$ with a sign-changing root $\alpha_i$ such that \[|\alpha_i-(\alpha_{i-1}+1)|<\eta\]

Consequently, \[|\alpha_k-(\alpha_0+k)|<k\eta\]

and thus, \begin{align*}|\alpha_k-x|&\leq|\alpha_k-(\alpha_0+k)|+|\alpha_0-y|+|y+k-x|\\&<(k+1)\eta+\frac{\epsilon}{2}<\epsilon\end{align*}

Throughout the proof, every digraph used was digon-free; thus, a digon-free digraph has a real dichromatic root within $\epsilon$ of $x$, and since $x$ and $\epsilon$ were arbitrary, it follows that real dichromatic roots of digon-free digraphs are dense in $\R$.
\end{proof}

\section{Dichromatic Uniqueness}
We say that a digraph $D$ is \textbf{dichromatically unique} if for every digraph $D'$, $P_D(x)=P_{D'}(x)$ implies $D\cong D'$.

\begin{figure}[h]
    \centering
    \resizebox{\textwidth}{!}{
        \input{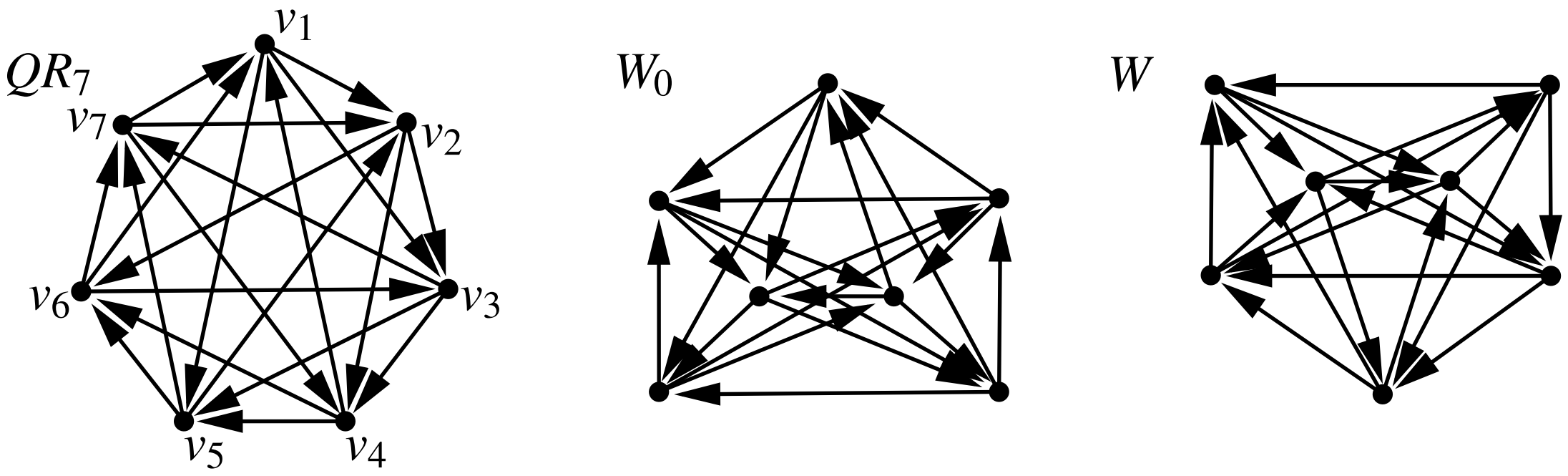}
    }
    \caption{The Tournaments $QR_7$, $W$, and $W_0$.}
    \label{fig:dic_unique}
\end{figure}

In \cite{DichromaticPoly}, the authors ask if the tournaments $QR_7$, $W$, and $W_0$ (see Figure \ref{fig:dic_unique}) are dichromatically unique, and we answer this affirmatively. We use the following theorem.
\begin{thm}\cite{ProjectGutenberg}\label{thm:tournament_triangle_bound}
    Let $c(T)$ denote the number of $3$-cycles in the tournament $T$ where $T$ has $n$ vertices. Then \[c(T)\leq \begin{cases}
        \frac{1}{24}(n^3-n) &\text{if $n$ is odd,}\\
        \frac{1}{24}(n^3-4n) &\text{if $n$ is even.}
    \end{cases}\]
    Equality holds throughout only for regular tournaments.
\end{thm}

\begin{thm}\label{thm:3_tournaments_unique}
The tournaments $QR_7$, $W$, and $W_0$ are dichromatically unique.
\end{thm}

\begin{proof}
$QR_7$, $W$, and $W_0$ are the only regular tournaments on $7$ vertices, so by \Cref{thm:tournament_triangle_bound}, we know that each of them has exactly $14$ directed triangles. \Cref{corr:third_coefficient} then implies that the third coefficient of their dichromatic polynomials equals $-14$. Using the code and database from \href{https://github.com/Katiekuehr/Dichromatic-Polynomial}{here} (see Appendix A), we can perform a brute-force search on all relevant tournaments to see that $QR_7$, $W$, and $W_0$ are dichromatically unique.
\end{proof}

\vspace{-0.4cm}

\section{Asymptotic Density of (Di)chromatic Polynomials}

Beyond determining the algebraic and combinatorial properties of chromatic polynomials, it is natural to ask how frequently polynomials with prescribed algebraic properties admit a coloring interpretation; in other words, what proportion of polynomials arise as the chromatic polynomial of a graph. This is the classical problem of characterizing chromatic polynomials, which has been studied by Wilf \cite{wilf1976} and subsequently placed in a general framework of recognition and characterization problems for graph polynomials by Kotek et. al. \cite{Kotek_Makowsky_Ravve_2018}. The density question considered in this section was motivated by earlier joint work of one of the authors on the proportion of irreducible polynomials in unicritically generated semigroups \cite{Hindes_Jacobs_Keller_Kim_Ye_Zhou_2025}. That work investigated how frequently a specific arithmetic property (irreducibility over $\Q$) occurred within a structured family of polynomials, which is in some ways analogous to the question we ask here about chromatic and dichromatic polynomials. In our case, the linear factorizations associated with chromatic polynomials suggests a family built from product factors $(x-j)$ for $j\in\Z_{\geq 0}$, which we can generalize slightly to include finitely many multiplicative and additive modifications. Within the families below, we show that the proportion of polynomials arising as dichromatic polynomials of strongly connected digraphs remains bounded away from both zero and one as the size of the digraphs (and thus degree of the polynomials) increases. Thus, realizable (di)chromatic polynomials constitute a substantial part of these algebraically defined families while the constraints imposed by a colouring interpretation exclude a substantial proportion as well. 

We begin with a lemma which allows us to consider dichromatic polynomials arising from strongly connected digraphs as a generalization of chromatic polynomials.

\begin{lem}\label{lem:strongly_connected_same_poly}
Let $G=(V,E)$ be a graph with at least one edge. Then there exists a strongly connected digraph $D=(V,A)$ (on the same vertex set) such that \[P_D(x)=\chi_G(x).\]
\end{lem}

\begin{proof}
In general, it is the case that \[P_{\bi G}(x)=\chi_G(x)\]

where $\bi G$ denotes the biorientation of $G$. In the case that $G$ is disconnected, then $\bi G$ is not strongly connected. So, choose some strongly connected component $C_0$ of $\bi G$ which contains the directed edge $(u,v)$. For each other strongly connected component $C_i$ of $\bi G$, choose a vertex $r_i\in C_i$ and add the edges $(u,r_i)$ and $(r_i,v)$. Clearly the resulting digraph is strongly connected, and, crucially, the new arcs do not change the valid colourings; any monochromatic obstruction still comes exactly from a monochromatic digon corresponding to an original graph edge. Therefore the dichromatic polynomial remains $\chi_G(x)$. This completes the proof.
\end{proof}

We now define the family of polynomials described briefly in the introduction to this section. For $\lambda>0$, define $m(n)=\lfloor\lambda\sqrt{n}\rfloor$, let $\m Q\subset\Z[x]$ be a finite set of nonzero monic polynomials containing $1$, and let $\m R\subset\Z[x]$ be a finite set of polynomials containing $0$ such that $\deg R<n$ for all $R\in\m R$. Now define \[\m F_n(\m Q,\m R,\lambda)=\left\{Q(x)\prod_{j=0}^{m(n)}(x-j)^{a_j}+R(x)\mid Q\in\m Q, R\in\m R, a_0\geq 1, a_j\geq0,\deg Q+\sum_{j=0}^{m(n)}a_j=n\right\}\]

Furthermore, let $S$ denote the class of all strongly connected digraphs and let $P(S)=\{P_D(x)\mid D\in S\}$ be the set of all dichromatic polynomials of strongly connected digraphs. Note that the choice of $m(n)=\lfloor\lambda\sqrt n\rfloor$ is solely to encourage a limiting proportion strictly between zero and one.

\begin{thm}[Density Theorem]\label{thm:asymptotic_density_theorem}
Fix $\m Q,\m R\subset\Z[x]$ according to the definitions above, fix $\lambda>0$, and let $\m F_n=\m F_n(\m Q,\m R,\lambda)$. Then,

\[0<\liminf_{n\to\infty}\frac{|\m F_n\cap P(S)|}{|\m F_n|}\leq\limsup_{n\to\infty}\frac{|\m F_n\cap P(S)|}{|\m F_n|}<1.\]

In other words, $P(S)$ has a positive, non-full asymptotic density in $\m F_n$.
\end{thm}

\begin{proof}
We first isolate the ``base layer" of $\m F_n$, which is simply the part of $\m F_n$ corresponding to $Q(x)=1$ and $R(x)=0$. Define \[\m S_n=\left\{\prod_{j=0}^{m(n)}(x-j)^{a_j}\mid a_0\geq 1, a_j\geq 0, \sum_{j=0}^{m(n)}a_j=n\right\}.\]

Clearly, $\m S_n\subseteq\m F_n$, and moreover, since $\m Q$ and $\m R$ are finite, the whole family $\m F_n$ is at most a bounded number of these ``layers." Specifically, for each pair $(Q,R\in\m Q\times\m R$, the number of possible exponent vectors is at most $|\m S_n|$. This tells us that \[|\m F_n|\leq |\m Q||\m R||\m S_n|\]

Also notice that (via stars and bars counting) \[|\m S_n|=\binom{n+m(n)-1}{m(n)}.\]

Call a polynomial in $\m S_n$ ``gap-free" if it has the form $\prod_{j=0}^{k}(x-j)^{a_j}$ for some $0\leq k\leq m(n)$ with $a_j\geq 1$ for all $j$. For a fixed $k$, the number of such polynomials is simply the number of positive solutions to $a_0+a_1+\cdots+ a_k=n$, which is $\binom{n-1}{k}$. Hence, the number of gap-free elements of $\m S_n$ is \[N_n=\sum_{k=0}^{m(n)}\binom{n-1}{k}.\]

We claim that all gap-free elements (except $x^n$) lie in $P(S)$. To see this, let \[f(x)=\prod_{j=0}^{k}(x-j)^{a_j}\text{ where } a_j\geq 1.\]

If $k\geq 1$, then we can construct a graph $G$ such that $\chi_G(x)=f(x)$ in the following way. Start with $K_{k+1}$, whose chromatic polynomial is $x(x-1)\cdots(x-k)$. Then, for each $j\in\{0,\dots,k\}$, add $a_j-1$ new vertices adjacent exactly to a fixed $j$-clique in the $K_{k+1}$ graph. Then each added type $j$ vertex contributes a factor of $x-j$. Call the resulting graph $G$. Then it is clear that \[\chi_G(x)=\prod_{j=0}^{k}(x-j)^{a_j}=f(x).\]

Furthermore, since $k\geq 1$, $G$ has at least one edge, so by \Cref{lem:strongly_connected_same_poly}, there exists a strongly connected digraph $D$ such that $P_D(x)=\chi_G(x)=f(x)$, so $f(x)\in P(S)$. The only excluded gap-free polynomial is $x^n$ (which corresponds to $k=0$), but this single polynomial is asymptotically negligible. Hence, we have that \[|\m S_n\cap P(S)|\geq N_n-1.\]

Combining this with previous results yields: \[\frac{|\m F_n\cap P(S)|}{|\m F_n|}\geq\frac{|\m S_n\cap P(S)|}{|\m Q||\m R||\m S_n|}\geq\frac{N_n-1}{|\m Q||\m R||\m S_n|}.\]

We now shift focus to computing the asymptotic ratio $N_n/|\m S_n|$ and demonstrate that the final term of $N_n$ dominates the entire sum asymptotically. For $1\leq k\leq m(n)$, we have that \[\frac{\binom{n-1}{k-1}}{\binom{n-1}{k}}=\frac{k}{n-k}\leq\frac{m(n)}{n-m(n)}\sim\frac{\lambda\sqrt{n}}{n}=\frac{\lambda}{\sqrt{n}}\overset{n\to\infty}{\to} 0\]

since $m(n)=\lfloor\lambda\sqrt{n}\rfloor$. In particular, if $\rho_n=\frac{m(n)}{n-m(n)}$, then for all $1\leq k\leq m(n)$, we have that \[\binom{n-1}{k-1}\leq\rho_n\binom{n-1}{k}.\]

In particular, this implies that for suitable $r$, \[\binom{n-1}{m(n)-r}\leq\rho_n^r\binom{n-1}{m(n)}.\]

Thus, we can bound $N_n$ as follows: \[N_n=\sum_{k=0}^{m(n)}\binom{n-1}{k}\leq\binom{n-1}{m(n)}\sum_{k=0}^{m(n)}\rho_n^k\leq\binom{n-1}{m(n)}\frac{1}{1-\rho_n}=(1+o(1))\binom{n-1}{m(n)}\]

since $\rho_n\to 0$ (throughout this proof, $o(1)$ is used to denote a quantity tending to $0$). This gives us that \[\binom{n-1}{m(n)}\leq N_n\leq (1+o(1))\binom{n-1}{m(n)}\]

which further implies that \[N_n=(1+o(1))\binom{n-1}{m(n)}.\]

Thus, we get that \[\frac{N_n}{|\m S_n|}=(1+o(1))\frac{\binom{n-1}{m(n)}}{\binom{n+m(n)-1}{m(n)}}=(1+o(1))\prod_{i=0}^{m(n)-1}\frac{n-i-1}{n+i}.\]

We can further determine the asymptotics of this product as follows. Taking the logarithm of the product yields \[\ln\left(\prod_{i=0}^{m(n)-1}\frac{n-i-1}{n+i}\right)=\sum_{i=0}^{m(n)-1}\ln\left(\frac{\frac{n-i-1}{n}}{\frac{n+i}{n}}\right)=\sum_{i=0}^{m(n)-1}\left[\ln\left(1-\frac{i+1}{n}\right)-\ln\left(1+\frac{i}{n}\right)\right]\]

\[=-\frac{1}{n}\sum_{i=0}^{m(n)-1}(2i+1)+O\left(\sum_{i=0}^{m(n)-1}\frac{i^2+1}{n^2}\right).\]

Here, the last equality holds since $m(n)=O(\sqrt{n})$, meaning that we can use the Taylor expansions $\ln(1-y)=-y+O(y^2)$ and $\ln(1+y)=y+O(y^2)$ where $y=i/n$ or $y=(i+1)/n$. Since $0\leq i\leq m(n)-1$ and $m(n)=\lfloor\lambda\sqrt{n}\rfloor$, we have that \[\frac{i}{n}\leq\frac{m(n)}{n}=O(n^{-1/2})\to 0.\]

So the Taylor expansions are valid uniformly over all $i$ in the sum. However, \[\sum_{i=0}^{m(n)-1}(2i+1)=(m(n))^2\]

and \[\sum_{i=0}^{m(n)-1}\frac{i^2+1}{n^2}=\frac{1}{n^2}\sum_{i=0}^{m(n)-1}(i^2+1)=\frac{1}{n^2}\left(m(n)+\sum_{i=0}^{m(n)-1}i^2\right)=O\left(\frac{(m(n))^3}{n^2}\right)=o(1).\]

Hence, \[\ln\left(\prod_{i=0}^{m(n)-1}\frac{n-i-1}{n+i}\right)=-\frac{(m(n))^2}{n}+o(1)\overset{n\to\infty}{\longrightarrow}-\lambda^2.\]

Exponentiating then implies that \[\frac{N_n}{|\m S_n|}\to e^{-\lambda^2}.\]

So, we conclude that \[\liminf_{n\to\infty}\frac{|\m F_n\cap P(S)|}{|\m F_n|}\geq\frac{e^{-\lambda^2}}{|\m Q||\m R|}>0.\]

This proves the lower-density bound. For the upper-density bound, we work with $\m S_n$ again. We say that a polynomial \[f(x)=\prod_{j=0}^{m(n)}(x-j)^{a_j}\in\m S_n\]

has a ``positive root gap" if there exist integers $0\leq i<j\leq m(n)$ such that $a_i=0$ and $a_j\geq 1$. Clearly, no polynomial with a positive root gap can be a dichromatic polynomial. This is because if $f(x)=P_D(x)$ for some digraph $D$, then $f(q)=0$ for all $q\in\{0,\dots,\vec\chi(D)\}$. However, this is not possible if $f$ has a positive root gap.

\vspace{10pt}

Thus, the only elements of $\m S_n$ without a positive root gap are the $N_n$ gap-free elements we counted above. Hence, at least $|\m S_n|-N_n$ elements of $\m S_n$ do not belong to $P(S)$. Now, since $\m S_n\subseteq\m F_n$, we have that \[|\m F_n\setminus P(S)|\geq |\m S_n|-N_n.\]

This gives us that \[\frac{|\m F_n\setminus P(S)|}{|\m F_n|}\geq\frac{|\m S_n|-N_n}{|\m Q||\m R||\m S_n|}.\]

Furthermore, we have that \[\frac{|\m S_n|-N_n}{|\m S_n|}\to 1-e^{-\lambda^2}.\]

This implies that \[\liminf_{n\to\infty}\frac{|\m F_n\setminus P(S)|}{|\m F_n|}\geq\frac{1-e^{-\lambda^2}}{|\m Q||\m R|}>0.\]

which is equivalent to saying that \[\limsup_{n\to\infty}\frac{|\m F_n\cap P(S)|}{|\m F_n|}\leq 1-\frac{1-e^{-\lambda^2}}{|\m Q||\m R|}<1.\]

This gives us the upper-density bound. Combining the lower and upper density bounds gives us the desired result, which is \[0<\liminf_{n\to\infty}\frac{|\m F_n\cap P(S)|}{|\m F_n|}\leq\limsup_{n\to\infty}\frac{|\m F_n\cap P(S)|}{|\m F_n|}<1.\]

In other words, the family of dichromatic polynomials coming from strongly connected digraphs has positive (but non-full) asymptotic density within the larger family $\m F_n$.
\end{proof}

The following corollary establishes one application of the above theorem.

\begin{corr}\label{density_application}
Fix $\lambda>0$ and let $m(n)=\lfloor\lambda\sqrt n\rfloor$. Define $\m H_n$ to be the set of distinct characteristic polynomials of nonempty central supersolvable hyperplane arrangements $\m A$ in $\R^n$ such that $1\leq\mathrm{rank}(\m A)<n$ and every exponent of the essentialization of $\m A$ is at most $m(n)$. Let $\m C$ denote the set of chromatic polynomials of (finite) simple graphs, and let $\m D_{\text{sc}}$ denote the set of dichromatic polynomials of (finite) strongly connected digraphs. Then, \[\lim_{n\to\infty}\frac{|\m H_n\cap\m C|}{|\m H_n|}=\lim_{n\to\infty}\frac{|\m H_n\cap\m D_{\text{sc}}|}{|\m H_n|}=e^{-\lambda^2}.\]

In particular, both limiting proportions lie strictly between zero and one, and this establishes a nontrivial density of (di)chromatic polynomials among characteristic polynomials of supersolvable arrangements.
\end{corr}

\begin{proof}
A supersolvable arrangement has a characteristic polynomial which factors into integer linear factors. Its positive roots are its nonzero exponents; a nonessential arrangement contributes factors of $x$, and a nonempty central arrangement contributes a factor of $x-1$ \cite{Stanley2007Arrangements}. Consequently, for $n\geq 2$ and $m(n)\geq 1$, \[\m H_n\subseteq\left\{\prod_{j=0}^{m(n)}(x-j)^{a_j}\mid a_0,a_1\geq 1, a_j\in\Z_{\geq 0}, \sum_{j=0}^{m(n)}a_j=n\right\}\]

The reverse inclusion holds as well. Indeed, take any polynomial $f(x)$ on the right, remove a factor of $(x-1)$, and list the remaining roots with multiplicity as $b_1,\dots,b_{n-1}$. In $\R^n$, with coordinates of $z,y_1,\dots,y_{n-1}$, take the hyperplanes $z=0$ and $y_i=kz$ for $1\leq i\leq n-1$ and $1\leq k\leq b_i$. For $b_i=0$, include no hyperplanes involving $y_i$. This arrangement is the cone over a product of one-dimensional arrangements, where the $i$-th factor consists of $b_i$ points. Hence its characteristic polynomial is \[(x-1)\prod_{i=1}^{n-1}(x-b_i)=f(x)\]

using the standard coning identity \cite{Stanley2007Arrangements}. This arrangement is also nonessential and supersolvable via Lemma 6.7 of \cite{Stanley2007Arrangements}. Hence, we have the following identification: \[\m H_n=\left\{\prod_{j=0}^{m(n)}(x-j)^{a_j}\mid a_0,a_1\geq 1, a_j\in\Z_{\geq 0}, \sum_{j=0}^{m(n)}a_j=n\right\}=\{f\in\m S_n\mid (x-1)\mid f\}\]

where we recall $\m S_n$ from the beginning of the proof of Theorem \ref{thm:asymptotic_density_theorem} above.

Consequently, $|\m H_n|=\binom{n+m(n)-2}{m(n)}$, and so, \[\frac{|\m H_n|}{|\m S_n|}=\frac{n-1}{n+m(n)-1}\longrightarrow 1\]

which implies that, asymptotically, $\m H_n$ contains all of $\m S_n$. So, for any set of polynomials $\m P$, \[\left|\frac{|\m H_n\cap\m P|}{|\m H_n|}-\frac{|\m S_n\cap\m P|}{|\m S_n|}\right|\leq\frac{|\m S_n\setminus\m H_n|}{|\m S_n|}=\frac{m(n)}{n+m(n)-1}\longrightarrow 0.\]

So, applying \Cref{thm:asymptotic_density_theorem} with $\m Q=\{1\}$ and $\m R=\{0\}$ immediately gives the desired result.
\end{proof}

\section{Open Conjectures/Problems}

\subsection{The Dichromatic Shameful Conjecture and Dichromatic Monotonicity}

For a graph $G$ of order $n$, the vertex-normalized function $\frac{\chi_G(q)}{q^n}$ is the probability that a uniformly random map $f: V(G) \rightarrow [q]$ of $V(G)$ is a colouring. In 1995, Bartels and Welsh conjectured that this probability cannot decrease when the number of available colours increases from $n-1$ to $n$, a conjecture which later became known as the \quotes{Shameful Conjecture} (\cite{Bartels_Welsh_1995}). This conjecture was proven in 2000 by Dong, together with the stronger result that $\frac{\chi_G(q)}{q^n}$ is monotonically increasing for integers $q\geq n-1$ (\cite{Dong_2000}). Since $\frac{P_D(q)}{q^n}$ likewise represents the probability that a uniformly random map $f: V(D) \rightarrow [q]$ of an order $n$ digraph $D$ contains no monochromatic directed cycle, it is natural to ask whether the analogous \quotes{Dichromatic Shameful Conjecture} holds for the dichromatic polynomial. As such, we propose the following: 

\begin{conj}[Dichromatic Shameful Conjecture]\label{conj:shameful_conj}
Let $D$ be a digraph of order $n$. Then \[\frac{P_D(n)}{n^n}\geq\frac{P_D(n-1)}{(n-1)^n}.\]
\end{conj}

Motivated by Dong's stronger result, it is also natural to ask whether the normalized dichromatic polynomial remains monotonically increasing beyond $n-1$ in general. As such, we propose the following stronger conjecture: 

\begin{conj}[Dichromatic Monotonicity Conjecture]\label{conj:monotonicity_conj}
Let $D$ be a digraph of order $n$. Then the function \[F_D(x)=\frac{P_D(x)}{x^n}\]

is monotonically increasing on $[n-1,\infty)$.
\end{conj}

With regard to \Cref{conj:monotonicity_conj}, it is clear that eventual dichromatic monotonicity holds automatically. Indeed, if $D$ is a digraph of order $n$, then there must exist $k\in\Z$ such that for all integers $x\geq k$, we have that $(x-1)^nP_D(x)\geq x^nP_D(x-1)$. Indeed, since both $(x-1)^nP_D(x)$ and $x^nP_D(x-1)$ are polynomials, we must have either $(x-1)^nP_D(x)\geq x^nP_D(x-1)$ or $(x-1)^nP_D(x)<x^nP_D(x-1)$ for sufficiently large $x$. If the latter is the case, then, adopting the notation from \Cref{conj:monotonicity_conj}, we can rewrite this inequality to be $F_D(x)<F_D(x-1)$. Since this inequality holds eventually, choose a large $N\in\N$ satisfying it. Note that by the definition of $F_D(x)$, $F_D(N)<F_D(N-1)\leq 1$. Since the sequence $F_D(N),F_D(N-1),\dots$ is decreasing, we have that \[\lim_{q\to\infty}F_D(q)\leq F_D(N)<1.\]

However, the limit exists and equals $1$ since $F_D$ is monic of degree $n$, which is a contradiction. Hence, the other inequality must hold eventually. The difficulty beyond eventual monotonicity is then showing that one may always take the uniform $k=n-1$ threshold.

\vspace{\baselineskip}

With regard to \Cref{conj:shameful_conj}, we note the following partial result. The classical Shameful Conjecture (now a theorem of Dong) implies that every graph $G$ of order $n$ (with $\chi_G(n-1)>0$) satisfies \[\frac{\chi_G(n)}{\chi_G(n-1)}>e,\] 

where $e=2.7182818284\dots$ is Euler's number. Furthermore, this bound is optimal although it is never attained \cite{Dong_2000}. Before the Shameful Conjecture was resolved, Seymour established the following nearly optimal bound: \[\frac{\chi_G(n)}{\chi_G(n-1)}\geq\frac{685}{252}=2.7182539682\dots\]

We show that the analogous bound holds for the dichromatic polynomial using Seymour's argument \cite{Seymour_1997a}. First, we recall a few of his definitions; we will use the term \textbf{potential $n$-colouring} of a graph $G$ to be \textit{any} map $f: V(G) \rightarrow [n]$, regardless of whether it admits monochromatic edges; the definition of a potential $n$-colouring of a digraph $D$ is analogous.

\begin{define}\label{define:consistent}
Let $G=(V,E)$ be a graph. We say a set $\m C$ of potential $n$-colourings of $V$ where $|V|=n$ is \textnormal{consistent} if it is closed under permutations of colours, and when we have a potential colouring of $G$ in $\m C$, and we choose a subset of one colour class and recolour it with some colour that does not yet appear, then this produces another potential colouring in $\m C$. More formally, $\m C$ is consistent if 

\begin{enumerate}
\item If $\phi\in\m C$ and $\pi$ is a permutation of $\{1,\dots,n\}$, then $\psi\in\m C$, where $\psi$ is defined by $\psi(v)=\pi(\phi(v))$ for all $v\in V$.
\item If $\phi\in\m C$, $1\leq i\leq n$, $X\subseteq\phi^{-1}(i)$, $1\leq j\leq n$, and $\phi^{-1}(j)=\emptyset$, then $\psi\in\m C$ where \[\psi(v)=\begin{cases}\phi(v)&\text{ if }v\notin X\\j&\text{ if }v\in X\end{cases}\]
\end{enumerate}
\end{define}

\begin{define}\label{define:expected_colours}
Let $\m C$ be a consistent set. Then define $r(\m C)$ to be the expected number of colours that appear in a randomly (uniformly) selected member of $\m C$.
\end{define}

We now show that the set of (acyclic) $n$-colourings of digraphs of order $n$ is a consistent set.

\begin{lem}\label{lem:acyclic_colourings_consistent}
Let $D=(V,A)$ be a digraph of order $n$, and let $C$ be the set of $n$-colourings of $D$. Then $C$ is consistent.
\end{lem}

\begin{proof}
Clearly, permuting the vertices of a colouring does not create a monochromatic directed cycle, so condition 1 is satisfied. Furthermore, since any induced subdigraph of an acyclic digraph is itself acyclic, recolouring some subset of a colour class to a completely new colour does not create a monochromatic cycle, so condition 2 is satisfied.
\end{proof}

We now obtain the following corollaries from Seymour: 

\begin{corr}\label{corr:expected_colours_bound}
Let $D=(V,A)$ be a digraph of order $n$, and let $C$ be the set of $n$-colourings of $D$. Then \[r(C)\geq n\left(1-\frac{252}{685}\right).\]
\end{corr}

\begin{proof}
Apply \Cref{lem:acyclic_colourings_consistent} and see the discussion on the top of page 196 in \cite{Seymour_1997a}. In particular, the argument depends only on the set of vertex colourings being consistent and does not reference edges, so does not use the underlying structure of the objects as graphs or digraphs.
\end{proof}

\begin{corr}\label{lem:expected_exact}
Let $D=(V,A)$ be a digraph of order $n$, and let $C$ be the set of $n$-colourings of $D$. Then \[r(C)=n\left(1-\frac{P_D(n-1)}{P_D(n)}\right)\]
\end{corr}

\begin{proof}
Apply \Cref{lem:acyclic_colourings_consistent} and see the discussion on the bottom of page 194 and the top of page 195 in \cite{Seymour_1997a}.
\end{proof}

Combining \Cref{lem:acyclic_colourings_consistent}, \Cref{corr:expected_colours_bound}, and \Cref{lem:expected_exact} allows us to prove the desired bound.

\begin{thm}\label{thm:shameful_partial_bound}
Let $D=(V,A)$ be a digraph of order $n$. Then \[P_D(n)\geq\frac{685}{252}P_D(n-1)\]
\end{thm}

\begin{proof}
Let $C$ be the set of $n$-colourings of $D$. By \Cref{lem:acyclic_colourings_consistent}, $C$ is consistent, so we can apply \Cref{corr:expected_colours_bound} and \Cref{lem:expected_exact} to get that \[n\left(1-\frac{P_D(n-1)}{P_D(n)}\right)\geq n\left(1-\frac{252}{685}\right)\implies\frac{P_D(n-1)}{P_D(n)}\leq\frac{252}{685}\]

This concludes the proof.
\end{proof} 

\section{AI Usage Declaration}
LLMs were used to assist with the code, and by one author to initially improve the writing for Section 4. No LLMs were used to generate mathematical arguments at any point.

\section*{Acknowledgments}

We would like to thank Sophie Spirkl for enlightening discussions.

\section*{Funding Declaration}

This research was funded by The Fields Institute as part of the 2026 Fields Undergraduate Summer Research Program.

We acknowledge the support of the Natural Sciences and Engineering Research Council of Canada (NSERC), [funding reference number RGPIN-2022-03093]. Cette recherche a \'et\'e financ\'ee par le Conseil de recherches en sciences naturelles et en g\'enie du Canada (CRSNG), [num\'ero de r\'ef\'erence RGPIN-2022-03093].

\printbibliography

\begin{appendices}
\section{Code}
The code for calculating the dichromatic polynomials of digraphs can be found in the GitHub Repository at \href{https://github.com/Katiekuehr/Dichromatic-Polynomial}{https://github.com/Katiekuehr/Dichromatic-Polynomial}.
\end{appendices}

\end{document}